\documentclass[envcountsect, envcountsame]{svmult}

\usepackage{mathptmx}       
\usepackage{helvet}         
\usepackage{courier}        
\usepackage{type1cm}        
\usepackage{makeidx}         
\usepackage{graphicx}        
\usepackage{multicol}        
\usepackage[bottom]{footmisc}

\makeindex             

\usepackage{amstext,amsfonts,amssymb,amscd,amsbsy,amsmath}
\usepackage{tikz}

\newcommand{\PP}{\mathbb{P}}

\newcommand{\RR}{\mathbb{R}}

\DeclareMathOperator{\Gr}{Gr}

\newcommand{\excise}[1]{}

\newcommand{\lra}[1]{\mathrel{\mathop{\longrightarrow}^{\mathrm{#1}}}} 

\smartqed

\begin{document}

\title*{The multidegree of the multi-image variety}
\author{Laura Escobar and Allen Knutson}
\institute{Laura Escobar \at 
University of Illinois at Urbana-Champaign, Department of Mathematics, Urbana, IL, 61801, USA \\ \email{lescobar@illinois.edu}
\and Allen Knutson \at Cornell University, Department of Mathematics, Ithaca, NY, 14850, USA, \\\email{allenk@math.cornell.edu}}
%
%
\maketitle

\abstract*{The multi-image variety is a subvariety of $\Gr(1,\PP^3)^n$ that models taking pictures with $n$ rational cameras. We compute its cohomology class in the cohomology of $\Gr(1,\PP^3)^n$, and from there its multidegree in the Pl\"ucker
embedding $(\mathbb{P}^5)^n$.}

\abstract{The multi-image variety is a subvariety of $\Gr(1,\PP^3)^n$ that models taking pictures with $n$ rational cameras. We compute its cohomology class in the cohomology of $\Gr(1,\PP^3)^n$, and from there its multidegree as a subvariety of $(\mathbb{P}^5)^n$ under the Pl\"ucker
embedding.}

\section{Introduction}

\emph{Multi-view geometry} studies the constraints imposed on a three-dimensional scene from various two-dimensional images of the scene.
Each image is produced by a camera.
\emph{Algebraic vision} is a recent field of mathematics 
in
which the techniques of algebraic geometry and optimization are used to formulate and solve problems in computer vision.
One of the main objects studied by this field is a \emph{multi-view variety.}
Roughly speaking, a multi-view variety parametrizes all the possible images that can be taken by a fixed collection of cameras. See \cite{AST,PST,THP} for more details on multi-view varieties.
See \cite{SRTGB} for a survey on various camera models.

The work \cite{PST} presents a new point of view to study the multi-view varieties.
A \emph{photographic camera} maps a point in the scene to a point in the image. \cite{PST} defines a \emph{geometric camera}, which maps a point in the scene not to a point, but to a viewing ray.
More precisely, a photographic camera is a map $\PP^3\dasharrow \PP^2$ or $\PP^3\dasharrow \PP^1\times\PP^1$, whereas a geometric camera is a map $\PP^3\dasharrow \Gr(1,\PP^3)$ such that the image of each point is a line containing the point.
The viewing ray corresponding to a point $p\in\PP^3$ is the line in $\Gr(1,\PP^3)$ this point gets mapped to.
 An important assumption is that light travels along rays which means that the image of a geometric camera is two-dimensional.

\begin{figure}[h]
\begin{center}
\includegraphics[scale=.45]{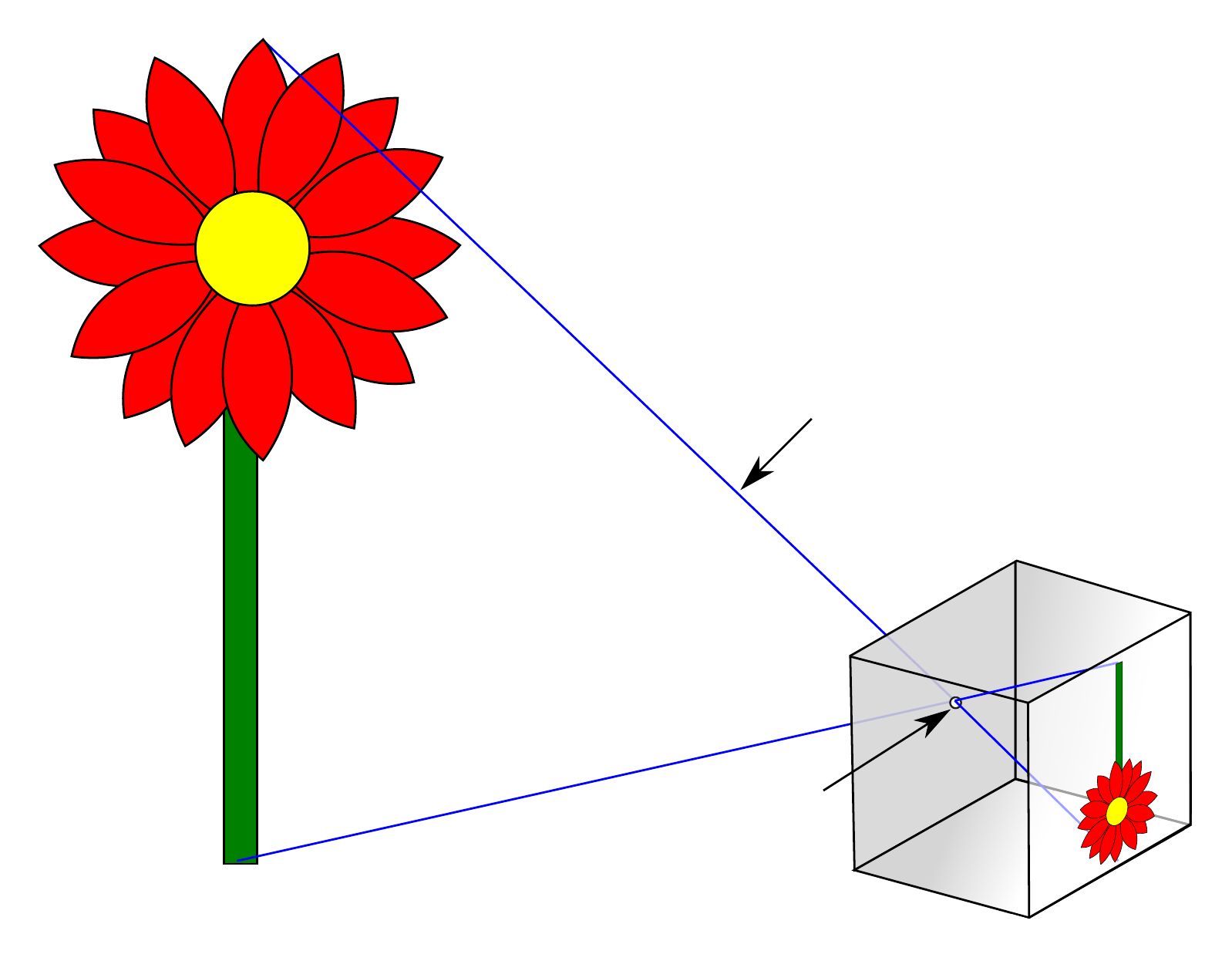}
\put(-165,160){$p$}
\put(-70,93){$\phi(p)$}
\put(-87,22){focal point}
\end{center}
\caption{A pinhole camera.}\label{fig:pinhole}
\end{figure}

We now illustrate these definitions using the example of \emph{pinhole cameras} or camer\ae\  obscur\ae, see Figure \ref{fig:pinhole}. 
As a geometric camera, a pinhole camera maps a point $p\in\PP^3$ to the line $\phi(p)$ connecting $p$ to the focal point.
This map is rational since it is undefined at the focal point.
As a photographic camera, a pinhole camera maps a point $p\in\PP^3$ to the intersection of $\phi(p)$ and the plane at the back of the camera.
Notice that, as a photographic camera, there are many pinhole cameras corresponding to a given focal point.
All of these cameras are equivalent up to a projective transformation.
The essential part of a pinhole camera is the mapping of the scene points to viewing rays, i.e. its modeling as a geometric camera.
 
A \emph{multi-image variety} is the Zariski closure of the image of a map 
\begin{align*}
\phi:\PP^3&\dasharrow \Gr(1,\PP^3)^n \\
p&\mapsto(\phi_1(p),\ldots,\phi_n(p)),\end{align*}
where each $\phi_i$ is a geometric camera.
The multi-view variety is a multi-image variety such that each $\phi_i$ is a pinhole camera.
Let $C_i$ be the closure of the image of $\phi_i$ inside the $i$-th $\Gr(1,\PP^3)$. 
Then under some assumptions, \cite[Theorem 5.1]{PST} shows that 
\begin{equation*}(C_1\times\cdots\times C_n)\cap V_n=\phi(\PP^3),
\end{equation*}
where $V_n$ is the \emph{concurrent lines variety} consisting of ordered $n$-tuples of lines in $\mathbb{P}^3$ that meet in a point $x$.
The concurrent lines and multi-image varieties are embedded into $(\PP^5)^n$ using the Pl\"ucker embedding.

The \emph{multidegree} of a variety embedded into a product of projective spaces is the polynomial whose coefficients give the numbers (when finite) of intersection points in the variety intersected with a product of general linear subspaces.
The present paper verifies the conjectured formula \cite[Equation (11)]{PST} for the multidegree of $V_n$ and computes the multidegree of the multi-image variety.
To do so, we describe $V_n$ using a projection of a partial flag variety and use Schubert calculus to compute the cohomology classes of $V_n$ and $(C_1\times\cdots\times C_n)\cap V_n$ in the cohomology ring of $\Gr(1,\PP^3)^n$. We then push forward these formul\ae\ into $(\PP^5)^n$ to obtain the multidegrees.

We now describe the organization of this paper.
In \S \ref{sec:miv} we define the main objects of study: the multi-image variety and the concurrent lines variety. We present the main theorem which computes the multidegrees of these objects.
The main tool to prove the main theorem is Schubert calculus.
In \S \ref{sec:schubertCalc} we give a brief introduction to Schubert calculus for $\Gr(1,\PP^3)$.
In \S \ref{sec:cohom} we compute the cohomology class of the multi-image variety and the concurrent lines variety in terms of the Schubert cycles in $\Gr(1,\PP^3)$.
We prove the main theorem by taking the pushforward of these equations to the cohomology ring of $\PP^5$.
In \S \ref{sec:Ktheory} we refine these results to a computation of the $K$-class for the concurrent lines variety.

\section{The multi-image variety}\label{sec:miv}

The \emph{Grassmannian} $\Gr(k,\PP^d)$ consists of $k$-dimensional planes inside $\PP^d$.
A \emph{congruence} is a two-dimensional family of lines in $\PP^3$, i.e. a surface $C$ in $\Gr(1,\mathbb{P}^3)$. The \emph{bidegree} $(\alpha,\beta)$ of a congruence $C$ is a pair of nonnegative integers such that the cohomology class of $C$ in $\Gr(1,\PP^3)$ has the form
\begin{equation}\label{eq:cong} [C] = \alpha \left[L:\ L\text{ contains a fixed point}\right]
  + \beta \left[L:\ L\text{ lies in a fixed plane}\right].
  \end{equation}
The first integer $\alpha$ is called the \emph{order} and counts the number of lines in $C$ that pass through a general point of $\PP^3$, and $\beta$ is called the \emph{class} and counts the number of lines in $C$ that lie in a general plane of $\PP^3$.
The \emph{focal locus} of $C$ consists of the points in $\PP^3$ that do not belong to $\alpha$ distinct lines of $C$.
\begin{example}A congruence with bidegree $(1,0)$ consists of all lines in $\Gr(1,\PP^3)$ that contain a fixed point. A geometric camera for such a congruence represents a pinhole camera where the fixed point is the focal point.
\end{example}

\begin{figure}[h]
\begin{center}
\includegraphics[scale=.45]{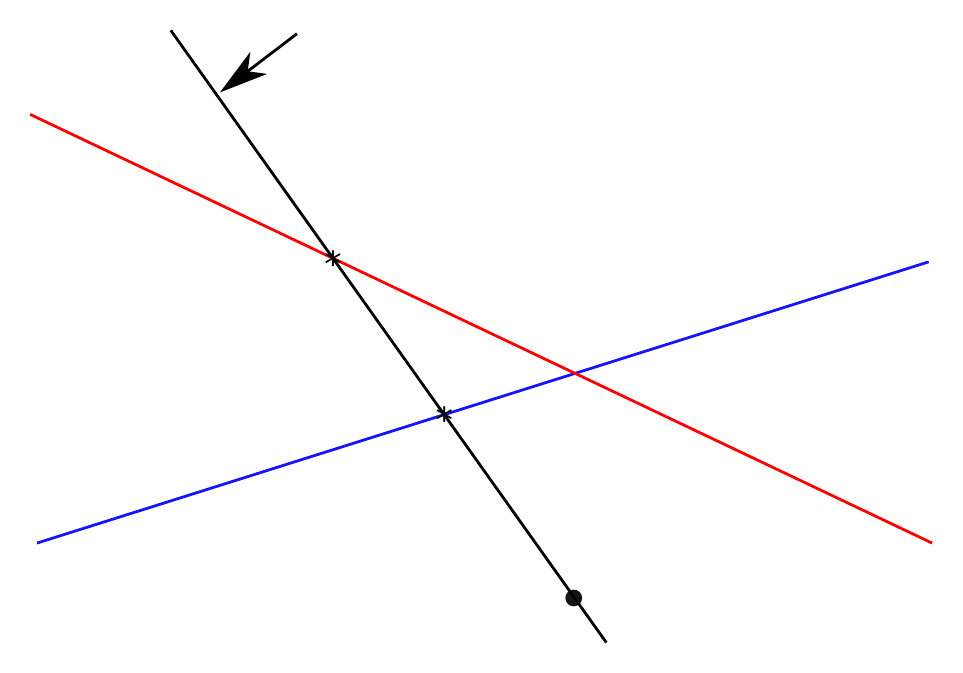}
\put(-60,8){$p$}
\put(-85,85){$\phi(p)$}
\put(-3,50){$L_1$}
\put(-3,15){$L_2$}
\end{center}
\caption{A two-slit camera.}\label{fig:2slit}
\end{figure}
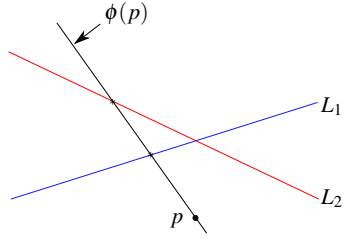

\begin{example}A two-slit camera assigns to $p\in\PP^3$ the unique line passing through $p$ and intersecting two fixed lines $L_1,L_2\in\Gr(1,\PP^3)$, see Figure \ref{fig:2slit}.
Its focal locus is $\{L_1,L_2\}$.
These cameras correspond to the congruences with bidegrees $(1,1)$.
\end{example}

The study of congruences started with \cite{Kum} which classified those of order one. They were studied by many mathematicians during the second half of the 19th century; see the book \cite{Jes} for some of these results.

\begin{remark} Congruences are studied the article \cite{KNT}. In particular, \cite[\S 2]{KNT} discusses congruences and their bidegrees. The bidegrees of curves in $\PP^1\times\PP^1$ are discussed in the article \cite{NFSY}.
\end{remark}

Consider a rational map 
\begin{align*}
\phi :  \mathbb{P}^3&\dasharrow C_1\times\cdots\times C_n\quad \subset \quad\Gr(1,\mathbb{P}^3)^n\\
x&\mapsto (\phi_1(x),\ldots,\phi_n(x))
\end{align*}
where the closure or $\phi_i(\PP^3)$ equals $C_i$, each $C_i$ is a congruence and $x\in \phi_i(x)$ for all $i$.
Each map $\phi_i$ is defined everywhere except on the focal locus of $C_i$.
In the language of \emph{algebraic vision} such a map means taking pictures with $n$ rational cameras where each $C_i$ is the $i$-th image plane.
An important assumption is that light travels along rays which means that $C_i$ is $2$-dimensional.

Let $C_1,\ldots,C_n$ be congruences of bidegree $(1,\beta_i)$.
Note that $\phi_i(x)$ is the unique line in $C_i$ passing through $x$.
The \emph{multi-image variety} of $(C_1,\ldots,C_n)$ is the Zariski closure of
$\phi(\PP^3)$. The \emph{concurrent lines variety} $V_n$ consists of ordered $n$-tuples of lines in $\mathbb{P}^3$ that meet in a point $x$.
By \cite[Theorem 5.1]{PST}, if the focal loci of the congruences are pairwise disjoint then the multi-image equals the intersection
$$(C_1\times\cdots\times C_n)\cap V_n\quad\subset\quad \Gr(1,\mathbb{P}^3)^n.$$

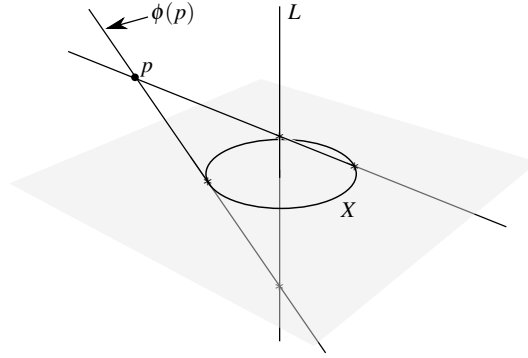
\begin{figure}[h]
\begin{center}
\includegraphics[scale=.6]{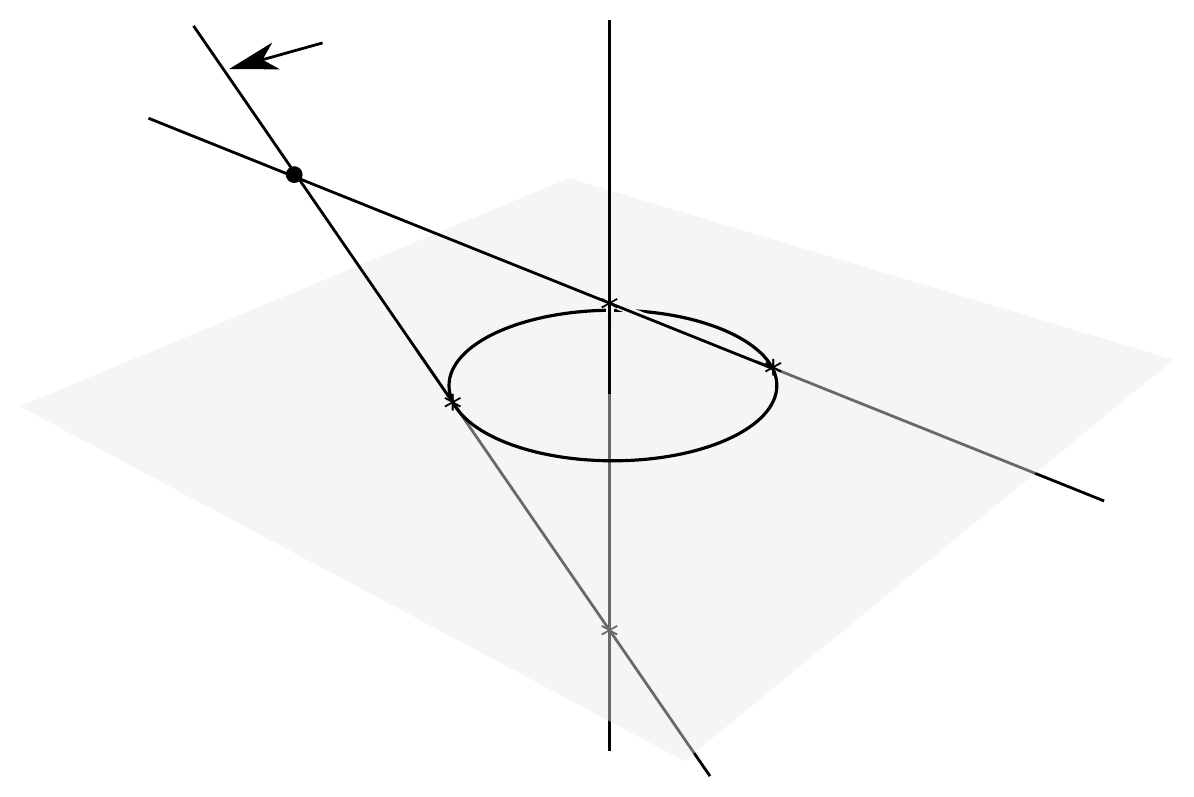}
\put(-98,130){$L$}
\put(-78,55){$X$}
\put(-154,110){$p$}
\put(-150,130){$\phi(p)$}
\end{center}
\caption{A non-central panoramic camera picks one of the two lines intersecting the circle $X$ and the line $L$.}\label{fig:pano}
\end{figure}

Most of the cameras studied in computer vision are associated with congruences of order $1$. However, cameras of higher order also appear in computer vision, see \cite{SRTGB} and \cite[\S 7]{PST}.
As an example of a camera associated to a congruence of bidegree $(2,2)$ let us discuss a \emph{non-central panoramic camera}, see Figure \ref{fig:pano}.
Consider a circle $X$ obtained by rotating a point about a vertical axis $L$.
There are two lines in $\Gr(1,\PP^3)$ passing through a general point $p\in\PP^3$ and intersecting both $L$ and $X$.
The congruence $C$ consisting of all lines intersecting both $X$ and $L$ has bidegree $(2,2)$.
A physical realization of a non-central panoramic camera consists of a sensor on the circle taking measurements pointing outwards.
This orientation of the sensor yields a map $\phi:\PP^3\dasharrow \Gr(1,\PP^3)$, i.e. it assigns only one line $\phi(p)$ to a point $p\in\PP^3$.

The concurrent lines and multi-image varieties are embedded into $(\PP^5)^n$ using the P\"ucker embedding
\begin{equation*}\Gr(1,\PP^3)\hookrightarrow \PP^5.
\end{equation*}
The \emph{multidegree} of a variety $X$ embedded into a product of projective spaces $\PP^{a_1}\times\cdots\times\PP^{a_n}$ is a homogeneous polynomial whose 
term $q\, z_1^{r_1}\cdots z_n^{r_n}$ indicates that there are $q$ intersection points when $X$ is intersected with the product $H_1\times\cdots\times H_n\subset \PP^{a_1}\times\cdots\times\PP^{a_n}$ of general linear subspaces, where $\dim(H_i)=r_i$.
The degree of this polynomial equals the codimension of $X$ in $\PP^{a_1}\times\cdots\times\PP^{a_n}$.
An equivalent definition of the multidegree of $X$ is its cohomology class in the cohomology ring of $\PP^{a_1}\times\cdots\times\PP^{a_n}$. The built-in command
\emph{multidegree} in the software \emph{Macaulay2} \cite{M2} computes the multidegree of $X$ from its defining ideal.
We refer to the book \cite[\S 8.5]{MS} for more details on multidegrees.

The main theorem of this paper computes the multidegree of the multi-image variety:
\begin{theorem}\label{thm:main}
The multidegree of the concurrent lines variety $V_n$ in $(\PP^5)^n$ equals
\begin{equation}\label{eq:condeg}(z_1z_2\cdots z_n)^3\left(4\sum_{(i,j) \atop i\neq j}z_i^{-2}z_j^{-1} + 8\sum_{\left\{ i,j,k\right\}}z_i^{-1}z_j^{-1}z_k^{-1} \right).\end{equation}
Let $(\alpha_i,\beta_i)$ be the bidegree of $C_i$ for $i=1,\ldots,n$. The multidegree of $(C_1\times\cdots\times C_n)\cap V_n$ in $(\PP^5)^n$ equals
\begin{align*}(\alpha_1\alpha_2\cdots\alpha_n)(z_1z_2\cdots z_n)^5
\Bigg(&\sum_{(i,j) \atop i\neq j}\frac{(\alpha_i+\beta_i)(\alpha_j+\beta_j)}{\alpha_i\alpha_j}z_i^{-2}z_j^{-1}+\\
 &\sum_{\left\{ i,j,k\right\}}\frac{(\alpha_i+\beta_i)(\alpha_j+\beta_j)(\alpha_k+\beta_k)}{\alpha_i\alpha_j\alpha_k}z_i^{-1}z_j^{-1}z_k^{-1}\Bigg),\end{align*}
where we distribute accordingly whenever $\alpha_i=0$.

\newpage
In particular, the multidegree of the multi-image variety of $(C_1,\ldots,C_n)$ where the bidegree of $C_i$ is $(1,\beta_i)$ equals
\begin{align*}(z_1z_2\cdots z_n)^5
\Bigg(&\sum_{(i,j) \atop i\neq j}(1+\beta_i)(1+\beta_j)z_i^{-2}z_j^{-1}+\\
 &\sum_{\left\{ i,j,k\right\}}(1+\beta_i)(1+\beta_j)(1+\beta_k)z_i^{-1}z_j^{-1}z_k^{-1}\Bigg).\end{align*}
\end{theorem}

\begin{remark} Ponce-Sturmfels-Trager \cite[equation (11)]{PST} had conjectured
  Formula \ref{eq:condeg} for the multidegree of $V_n$, based on
  experimental evidence from Macaulay2.
\end{remark}
\begin{remark} Aholt-Sturmfels-Thomas \cite[Corollary 3.5]{AST} gives an equation for the multidegree of the multi-image variety for the congruences with bidegree $(1,0)$. This equation coincides with the equation of Theorem \ref{thm:main}.
\end{remark}

To prove this theorem, we first use Schubert calculus to obtain the
cohomology classes of $V_n$ and the multi-image variety
in the cohomology ring of $\Gr(1,\mathbb{P}^3)^n$; see Theorems
\ref{thm:cohomVn} and \ref{thm:cohomMIV}. Using this formula we then describe the
multidegrees of these varieties and prove Theorem \ref{thm:main}.
In the next section we introduce our notation for Schubert varieties
and review the part of Schubert calculus that we need.

\section{Schubert varieties}\label{sec:schubertCalc}

In this section we review Schubert varieties in $\Gr(k,\PP^d)$ keeping $\Gr(1,\PP^3)$ as the main example. 
For further details we recommend the book \cite{Ful}.
\excise{Any point in $\Gr(k,\mathbb{P}^d)$ can be represented by a $(k+1)\times (d+1)$ matrix of rank $k+1$ in row echelon form. 
Let $C_{\left\{j_1,\ldots,j_{k+1}\right\}}$ be the \emph{Schubert cell} consisting of the $k$-dimensional planes that can be represented by matrices in row echelon form that have the pivots in positions $(1,j_1),(2,j_2),\ldots,(k+1,j_{k+1})$.
The Grassmannian has a cell decomposition
$$\Gr(k,\PP^n)=\bigsqcup_{\emph{j}\in \binom{[n+1]}{k+1}}C_\emph{j}.$$}
Fix a coordinate system for $\mathbb{P}^d$. 
Let $\mathbb{P}^{(i,i+1,\ldots,j)}$ denote the coordinate subspace of $\mathbb{P}^n$ spanned by the coordinates $(i,i+1,\ldots,j)$ and consider the \emph{standard flag}
\begin{equation*} {\sf E}_\bullet:=\PP^{(0)}\subset \PP^{(0,1)}\subset\cdots\subset \PP^{(0,1,\ldots,d)}=\mathbb{P}^d.
\end{equation*}
\excise{All the points in $C_\emph{j}$ intersect the planes in ${\sf F}_\bullet$ with dimensions prescribed by $\emph{j}$, i.e.
\begin{equation*}C_\emph{j}=\left\{ p\in\Gr(k,\PP^n) : \dim(p\cap \langle e_0,e_q,\ldots,e_{j_l-1} \rangle)= l \text{ for } l=1,\ldots,k+1 \right\}.
\end{equation*}}
The \emph{Schubert variety} in $\Gr(k,\PP^d)$ corresponding to the subset $\{j_1,\ldots,j_{k+1}\}\subset\{1,2,\ldots,d+1\}$ is
\begin{equation}\label{eq:schuvar}
X_{\{j_1,\ldots,j_{k+1}\}}:=\left\{ p\in\Gr(k,\PP^d) : \dim(p\cap \PP^{(0,1,\ldots,j_l-1)})\geq l-1 \text{ for } l=1,\ldots,k+1 \right\}.
\end{equation}

\begin{example} The $\binom{d+1}{k+1} = {4\choose 2}$ Schubert varieties in $\Gr(1,\PP^3)$ are
\begin{align*}
X_{\left\{1,2\right\}}&=\{\PP^{(0,1)} \},\\
X_{\left\{1,3\right\}}&=\{L\in \Gr(1,\PP^3) : \PP^{(0)}\subset L \subset \PP^{(0,1,2)}  \},\\
X_{\left\{1,4\right\}}&=\{L\in \Gr(1,\PP^3) : \PP^{(0)}\subset L  \},\\
\end{align*}
\begin{align*}
X_{\left\{2,3\right\}}&=\{L\in \Gr(1,\PP^3) : L\subset \PP^{(0,1,2)}  \},\\ 
X_{\left\{2,4\right\}}&=\{L\in \Gr(1,\PP^3) : \dim(L\cap \PP^{(0,1)}) \geq 1  \},\text{ and}\\ 
 X_{\left\{3,4\right\}}&=\Gr(1,\PP^3).
\end{align*}
\end{example}

Consider the Schubert cells defined by replacing $\geq$ by $=$ in Equation (\ref{eq:schuvar}).
Since $\Gr(1,\PP^3)$ is a disjoint union of the Schubert cells,
and they are each contractible, 
the classes $[X_J]$ for $J\subset\{1,2,\ldots,d+1\}$ form a basis for the cohomology ring of $\Gr(k,\PP^d)$. The ring operation is the cup product
\begin{equation*}
[X_I]\smallsmile[X_J]:=[X_I({\sf E}^{\text{op}}_\bullet)\cap X_J],
\end{equation*}
where $X_I({\sf E}^{\text{op}}_\bullet)$ is defined by using $\PP^{(d-j_l+1,d-j_l+2,\ldots,d)}$ instead of $\PP^{(0,1,\ldots,j_l-1)}$ in Equation \ref{eq:schuvar}.
(Unlike $X_I$, this variety intersects $X_J$ transversely, while having the same
cohomology class as $X_I$.)
One also gets a basis for the cohomology ring of products of Grassmannians via the K\"unneth isomorphism.

\begin{example}\label{ex:easycup} Note that $[X_J]\smallsmile[X_{\{1,2\}}]=0$ for any $J\neq\{3,4\}$ since $\PP^{(2,3)}\notin X_J$. Any line $L$ containing $\PP^{(0)}$ is not contained in $\PP^{(1,2,3)}$ and therefore $[X_{\{1,4\}}]\smallsmile[X_{\{2,3\}}]=0$. We have that $[X_{\{1,4\}}]\smallsmile[X_{\{1,4\}}]=[X_{\{1,2\}}]$ since there is a unique line containing the points $\PP^{(0)}$ and $\PP^{(3)}$. 
On general Grassmannians, computing $[X_I]\smallsmile[X_J]$ can be done using the ``Pieri rule'' for special classes, see Equation (\ref{eq:Pieri}), and the ``Littlewood-Richardson rule'' \cite{LR} for arbitrary dimensions.
\end{example}

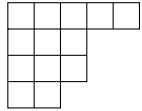
\begin{figure}[h]\centering
\begin{tikzpicture}[scale=0.7 ]
\draw (1,2)--(1,0)--(0,0)--(0,2)--(2.5,2)--(2.5,1.5)--(0,1.5);
\draw (2,2)--(2,1.5);
\draw (1.5,2)--(1.5,.5)--(0,.5);
\draw (.5,2)--(.5,0);
\draw (0,1)--(1.5,1);
\end{tikzpicture}
\caption{The Young diagram of the partition $\lambda=(5,3,3,2)$}\label{fig}
\end{figure}

Schubert varieties stratify $\Gr(k,\mathbb{P}^d)$.
The poset of their inclusions is most easily described 
when Schubert varieties are indexed using partitions.
A \emph{partition} of $n$ into $k+1$ parts is a list $\lambda=(\lambda_1\geq\lambda_2\geq\cdots\geq \lambda_{k+1} > 0)$ such that $n=\sum_i\lambda_i$. There is a bijection between $(k+1)$-subsets of $[d+1]$ and partitions $\lambda$ with at most $k+1$ parts such that $\lambda_1\leq d-k$ which is given by
\begin{equation*}
\emph{j}=\{j_1,\ldots,j_{k+1}\} 
\  \ \leftrightarrow\  \
\lambda=(d-k+1-j_1,\ d-k+2-j_2, \ldots,\ d-j_{k},d+1-j_{k+1}).
\end{equation*}
Partitions can be visualized in the following way. Given $\lambda_1\geq\lambda_2\geq \ldots\geq\lambda_{k+1}$, we draw a figure made up of squares sharing edges that has $\lambda_1$ squares in the first row, $\lambda_2$ squares in the second row starting out right below the beginning of the first row, and so on. 
This figure is called a \emph{Young diagram}. 
See Figure \ref{fig} for an example. 
We say that $\lambda\subset \lambda'$ if $\lambda_i\leq \lambda'_i$ for all $i$, i.e. if the diagram of $\lambda$ lies inside the one for $\lambda'$.

From now on we will index Schubert cells and varieties by partitions. With this indexing set we have the following facts:
\begin{itemize} 
\item $\dim(X_\lambda) = (k+1)(d-k)-\sum_i\lambda_i$,
\item $X_\lambda \supset X_\mu$ iff $\mu\supset\lambda$, and
\item \emph{The Pieri rule:} suppose that $\lambda=(\lambda_1,0,\ldots,0)$ and $\mu=(\mu_1,\mu_2,\ldots,\mu_{k+1})$ is any partition. Then
\begin{equation}\label{eq:Pieri}
[X_\lambda]\smallsmile[X_\mu]=\sum [X_\nu],
\end{equation}
where the sum is over all partitions $\nu=(\nu_1,\nu_2,\ldots,\nu_{k+1})$ such that $\nu_1\leq d-k$, $\mu_i\leq \nu_i\leq \mu_{i-1}$, and $\sum\nu_i=\sum(\lambda_i+\mu_i)$.
\end{itemize}
\begin{example}\label{ex:piericup} The Schubert varieties in $\Gr(1,\mathbb{P}^3)$ are ordered by containment in the poset in Figure \ref{fig:poset}.
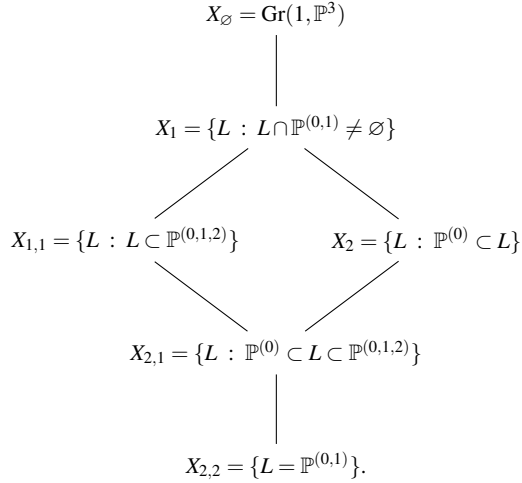
\begin{figure}[h]\begin{center}
\begin{tikzpicture}
\node (22) at (0,0) {$X_{2,2}=\{L = \PP^{(0,1)}\}$.};
\node (21) at (0,1.5) {$X_{2,1}=\{L\ :\ \PP^{(0)}\subset L\subset \PP^{(0,1,2)}\}$};
\node (11) at (-2,3) {$X_{1,1}=\{L\ :\ L\subset \PP^{(0,1,2)}\}$};
\node (2) at (2,3) {$X_{2}=\{L\ :\ \PP^{(0)}\subset L\}$};
\node (1) at (0,4.5) {$X_1=\{L\ :\ L\cap \PP^{(0,1)}\neq\varnothing\}$};
\node (0) at (0,6) {$X_{\varnothing}=\Gr(1,\PP^3)$};
\draw (22)--(21);
\draw (21)--(11);
\draw (21)--(2);
\draw (11)--(1);
\draw (2)--(1);
\draw (1)--(0);
\end{tikzpicture}\end{center}\caption{The containment poset of the Schubert varieties in $\Gr(1,\PP^3)$.}\label{fig:poset}\end{figure}
This poset is ranked by the dimensions of the Schubert varieties.
By the Pieri rule, we have that $[X_1]\smallsmile[X_2]=[X_{2,1}]$ and $[X_1]\smallsmile[X_{1,1}]=[X_{2,1}]$.
\end{example}

\begin{remark} We can rewrite Equation \ref{eq:cong} for the cohomology class of a congruence as
\begin{equation*} [C] =\alpha[X_2]+\beta[X_{1,1}]. 
  \end{equation*}
\end{remark}

\begin{remark} In the article \cite[\S 6]{KNT}, Schubert varieties and the intersection theory of $\Gr(1,\PP^3)$ are also discussed.
\end{remark}

\section{Computing the multidegrees}\label{sec:cohom}

Let $M\subset {\Gr}(0,\mathbb{P}^3)\times \Gr(1,\mathbb{P}^3)$ be the partial
flag manifold consisting of pairs $(P,L)$ where $P$ is a point in
$\mathbb{P}^3$ and $L$ is a line through $P$. Such a pair is called a
\emph{(partial) flag} in $\mathbb{P}^3$. Define $M^n_\Delta$ to be
the subvariety of $M^n$ consisting of lists of $n$ flags such that the point
$P$ is the same in all of them. Consider the diagonal
\begin{equation}
(\mathbb{P}^3)^n_\Delta:=\{(P_1,\ldots,P_n) : P_1=\cdots=P_n\};
\end{equation}
then $M^n_\Delta$ is the preimage of $(\mathbb{P}^3)^n_\Delta$ under the projection $M^n\rightarrow (\PP^3)^n$ induced by 
$$  M\longrightarrow \mathbb{P}^3, \qquad
(P,L)\longmapsto P. $$
For example, $M^2_\Delta \subseteq M^2$
consists of pairs of flags of the form $((P,L_1),(P,L_2))$.  
The following is straightforward:
\begin{proposition}The concurrent lines variety $V_n$ is the image of $M^n_\Delta$ under the projection $\displaystyle{M^n\lra{p} \Gr(1,\mathbb{P}^3)^n}$ induced by
$$M \hookrightarrow \Gr(0,\mathbb{P}^3)\times \Gr(1,\mathbb{P}^3)\lra{}  \Gr(1,\mathbb{P}^3).$$
\end{proposition}
The intersection $\left( M^{i-1} \times M^2_\Delta \times M^{n-i-1} \right)\cap\left( M^{i} \times M^2_\Delta \times M^{n-i-2} \right)$ consists of the lists of flags $((P_1,L_1),\ldots,(P_n,L_n))$ such that $P_i=P_{i+1}=P_{i+2}$, and this intersection is transverse.
We can write
$M^n_\Delta$ as the transverse intersection
\begin{equation}\label{eq:in} M^n_\Delta = \bigcap_{i=1}^{n-1} \left( M^{i-1} \times M^2_\Delta \times M^{n-i-1} \right).\end{equation}
From this description we deduce the following Theorems which give the cohomology classes of $V_n$ and the multi-image variety.
\begin{theorem}\label{thm:cohomVn}
  The class $[V_n]$ of the concurrent lines variety in the cohomology
  ring of $\Gr(1,\mathbb{P}^3)^n$ is
\begin{equation}\label{eq:equicoho}
 [V_n] =
  \sum_{0 = v_0 \leq v_1 \leq \ldots \leq v_n = 3
  \atop v_{j+1}-v_j<3} \ 
  \bigotimes_{i=0}^{n-1}
  \left[L: L\cap \mathbb{P}^{(v_{i},v_{i}+1,\ldots,v_{i+1})} \neq \varnothing \right].
\end{equation}
 Moreover, $[V_n]$ can be written as
\begin{eqnarray}\label{eq:multivn}
  [V_n] &=&
  \sum_{0 = v_0 \leq v_1 \leq \ldots \leq v_n = 3
  \atop v_{j+1}-v_j<3} \ 
  \bigotimes_{i=1}^n
  \begin{cases}
    [X_2] &\text{if }v_i-v_{i-1} = 0 \\
    [X_1] &\text{if }v_i-v_{i-1} = 1 \\
    [\Gr(1,\PP^3)] &\text{if }v_i-v_{i-1} = 2 
  \end{cases}.
\end{eqnarray}
\end{theorem}

\begin{remark}Equation (\ref{eq:equicoho}) is more natural than Equation (\ref{eq:multivn}),
  in that (\ref{eq:equicoho}) is correct in $T$-equivariant cohomology,
  whereas (\ref{eq:multivn}) is only correct in ordinary cohomology.
\end{remark}

\newpage
\begin{example} Consider $n=3$. In $T$-equivariant cohomology we have that
\begin{align*}[V_n]=&\left[L: L\cap \mathbb{P}^{(0)} \neq \varnothing \right]\otimes\left[L: L\cap \mathbb{P}^{(0,1)} \neq \varnothing \right]\otimes\left[L: L\cap \mathbb{P}^{(1,2,3)} \neq \varnothing \right]+\\
&\left[L: L\cap \mathbb{P}^{(0)} \neq \varnothing \right]\otimes\left[L: L\cap \mathbb{P}^{(0,1,2)} \neq \varnothing \right]\otimes\left[L: L\cap \mathbb{P}^{(2,3)} \neq \varnothing \right]+\\
&\left[L: L\cap \mathbb{P}^{(0,1)} \neq \varnothing \right]\otimes\left[L: L\cap \mathbb{P}^{(1)} \neq \varnothing \right]\otimes\left[L: L\cap \mathbb{P}^{(1,2,3)} \neq \varnothing \right]+\\
&(4\text{ other terms}).
\end{align*}
In ordinary cohomology we have that 
\begin{align*}[V_n]=&[X_2]\otimes[X_1]\otimes[\Gr(1,\PP^3)]\ +\ [X_2]\otimes[\Gr(1,\PP^3)]\otimes[X_1]\ +\ \\
&[X_1]\otimes[X_2]\otimes[\Gr(1,\PP^3)]\ +\ [X_1]\otimes[X_1]\otimes[X_1]\ +\ [X_1]\otimes[\Gr(1,\PP^3)]\otimes[X_2]\ +\ \\
&[\Gr(1,\PP^3)]\otimes[X_2]\otimes[X_1]\ +\ [\Gr(1,\PP^3)]\otimes[X_1]\otimes[X_2].\\
\end{align*}
\end{example}

\begin{proof}
From equation (\ref{eq:in}) we have that 
$$[M^n_\Delta]=\prod_{i=1}^{n-1}
\left(1\otimes\cdots\otimes 1\otimes[M^2_\Delta]\otimes1\otimes\cdots\otimes1\right)$$
in the cohomology of $(\Gr(0,\mathbb{P}^3)\times \Gr(1,\mathbb{P}^3))^n$. 
Identifying $H^*((\mathbb{P}^3)^2) \cong
H^*(\mathbb{P}^3) \otimes H^*(\mathbb{P}^3)$
using the K\"unneth isomorphism, we have
\begin{eqnarray*}
  \left[ (\mathbb{P}^3)^2_\Delta\right]
  &=& \sum_{0\leq v \leq 3} \left[\mathbb{P}^{(0,1,\ldots v)}\right] \otimes \left[\mathbb{P}^{(v,v+1,\ldots 3)}\right]
\end{eqnarray*}
and hence
\begin{eqnarray*}
  \left[ (\PP^3)^n_\Delta \right] &=&
  \sum_{0 = v_0 \leq v_1 \leq \ldots \leq v_n = 3}
  \left[\mathbb{P}^{(v_0,\ldots, v_1)}\right]
   \otimes \cdots \otimes 
    \left[\mathbb{P}^{(v_{n-1},\ldots, v_n)}\right].
\end{eqnarray*}
Pulling this back to $M^n$, we get essentially the same formula
\begin{equation*}
  \left[ M^n_\Delta \right] =
  \sum_{0 = v_0 \leq v_1 \leq \ldots  \leq v_n = 3}
  \left[(P,L): P\in \mathbb{P}^{(v_0,\ldots, v_1)}\right] \otimes \cdots\otimes 
  \left[(P,L): P\in \mathbb{P}^{(v_{n-1},\ldots, v_n)}\right].
\end{equation*}

Pushing forward this class under the projection map $\displaystyle{M^n\lra{p} \Gr(1,\mathbb{P}^3)^n}$, we get $p_*[ M^n_\Delta] = [V_n]$.
Not all terms of the $[M^n_\Delta]$ formula above 
survive when we project $M^n_\Delta$ to $\Gr(1,\mathbb{P}^3)^n$ by forgetting $P$:
\begin{itemize}
\item The image of $\{(P,L):\ P\in \mathbb{P}^{(0,1,2,3)}\}$ under $\mathrm{p}$ equals $\Gr(1,\PP^3)$. Since for $v_i-v_{i-1}=3$ the dimension of $\{(P,L):\ P\in \mathbb{P}^{(0,1,2,3)}\}$ drops under $\mathrm{p}$, any factor $\left[(P,L):\ P\in \mathbb{P}^{(0,1,2,3)} \right]$ pushes down to $0$.
\item The other terms
$\left[(P,L):\ P\in \mathbb{P}^{(i,i+1,\ldots,j)} \right]$ push down to
$\left[L:\ L\cap  \mathbb{P}^{(i,i+1,\ldots,j)}  \neq \varnothing \right]$. 
\end{itemize}
It follows that
\begin{eqnarray*}
  [V_n] &=&
  \sum_{0 = v_0 \leq v_1 \leq \ldots \leq v_n = 3\atop v_{j+1}-v_j<3}
  \left[L:\ L\cap \mathbb{P}^{(v_0,\ldots, v_1)} \neq \varnothing \right]
  \otimes \cdots \otimes
  \left[L:\ L\cap \mathbb{P}^{(v_{n-1},\ldots, v_n)}\neq \varnothing \right].
\end{eqnarray*}
Note that for any $L\in\Gr(1,\mathbb{P}^3)$ we have that $L\cap \mathbb{P}^{(0,1,2)}\neq \varnothing$ and $L\cap\mathbb{P}^{(1,2,3)}\neq \varnothing$. Therefore the terms corresponding to  $ \mathbb{P}^{(0,1,2)} $ and $\mathbb{P}^{(1,2,3)} $ push down to $[\Gr(1,\mathbb{P}^3)]$.
\qed
\end{proof}

\begin{theorem}\label{thm:cohomMIV}
For $i=1,\ldots,n$, let $C_i \subset \Gr(1,\PP^3)$ be a general congruence with class
$\alpha_i [X_{2}] + \beta_i [X_{1,1}] \in H^4(\Gr(1,\PP^3))$. Then
\begin{eqnarray*}
[(C_1\times\cdots\times C_n)\cap V_n] 
  &=& 
  \sum_{0 = v_0 \leq v_1 \leq \ldots \leq v_n = 3
\atop v_{j+1}-v_j<3} \ 
  \bigotimes_{i=1}^n
  \begin{cases}
    \alpha_i [X_{2,2}] &\text{if }v_{i} = v_{i-1} \\
    (\alpha_i+\beta_i) [X_{2,1}] &\text{if }v_{i}= v_{i-1} + 1  \\
    \alpha_i [X_2] +\beta_i [X_{1,1}] &\text{if }v_{i}= v_{i-1} + 2
  \end{cases}.
\end{eqnarray*}
\end{theorem}

\begin{example} Let $C_1,C_2,C_3$ be three congruences with bidegrees $(\alpha_i,\beta_i)_{i=1,2,3}$. Then 
\begin{align*}[(C_1\times C_2\times C_3)\cap V_3] 
  =\ &\alpha_1[X_{2,2}]\otimes(\alpha_2\ +\ \beta_2)[X_{2,1}]\otimes(\alpha_3[X_2]\ +\ \beta_3[X_{1,1}])\ +\ \\
  &\alpha_1[X_{2,2}]\otimes(\alpha_2[X_2]\ +\ \beta_2[X_{1,1}])\otimes(\alpha_3\ +\ \beta_3)[X_{2,1}]\ +\ \\
  &(\alpha_1\ +\ \beta_1)[X_{2,2}]\otimes\alpha_2[X_{2,2}]\otimes(\alpha_3[X_2]\ +\ \beta_3[X_{1,1}])\ +\ \\
  &(4\text{ other similar terms}).
\end{align*}

\end{example}

\begin{proof}
Let $C_i \subset \Gr(1,\PP^3)$ be a surface of class 
$\alpha_i [X_{1,1}] + \beta_i [X_2] \in H^4(\Gr(1,\PP^3))$.
By Theorem \ref{thm:cohomVn},
\begin{align*}
    [(C_1\times\cdots\times C_n)&\cap V_n] \\ &=[V_n ] \smallsmile ([C_1]\otimes \cdots \otimes [C_n]) \\
  &=\sum_{0 = v_0 \leq v_1 \leq \ldots \leq v_n = 3
\atop v_{j+1}-v_j<3} \ 
  \bigotimes_{i=1}^n
  \left[L: L\cap \PP^{(v_{i-1},\ldots, v_i)} \neq \varnothing \right]
  \smallsmile
  (\alpha_i [X_{2}] + \beta_i [X_{1,1}]).
\end{align*}
Using the computations from Examples \ref{ex:easycup} and \ref{ex:piericup} we have that:
\begin{itemize}
\item When $v_{i} = v_{i-1}$, 
then $ \left[L:\ L\supset \PP^{(v_{i})}\right]=[X_2]$, 
\\ and the factor is $[X_2] \smallsmile
    (\alpha_i [X_{2}] + \beta_i [X_{1,1}])=\alpha_i [X_{2,2}].$
\item When $v_{i} = v_{i-1} + 1$, then
$ \left[L:\ L\cap \PP^{(v_{i-1}, v_{i})}\neq\varnothing\right]=[X_1]$, \\ 
and the factor is $[X_1]\smallsmile
    (\alpha_i [X_{2}] + \beta_i [X_{1,1}]) = (\alpha_i+\beta_i) [X_{2,1}] $.
\item When $v_{i} = v_{i-1}+2$, then
$ \left[L:\ L\cap \PP^{(v_{i-1},v_{i-1}+1, v_{i})}\neq\varnothing\right] =[\Gr(1,\mathbb{P}^3)]= 1$,
\\ giving $1 \smallsmile
    (\alpha_i [X_{2}] + \beta_i [X_{1,1}]) = 
    \alpha_i [X_{2}] + \beta_i [X_{1,1}]$.
\end{itemize}
\newpage
The result is
\begin{eqnarray*}
      [(C_1\times\cdots\times C_n)\cap V_n] 
  &=& 
  \sum_{0 = v_0 \leq v_1 \leq \ldots \leq v_n = 3
\atop v_{j+1}-v_j<3} \ 
  \bigotimes_{i=1}^n
  \begin{cases}
    \alpha_i [X_{2,2}] &\text{if }v_{i} = v_{i-1} \\
    (\alpha_i+\beta_i) [X_{2,1}] &\text{if }v_{i}= v_{i-1} + 1  \\
   \alpha_i [X_{2}] + \beta_i [X_{1,1}]&\text{if }v_{i}= v_{i-1} + 2.  
  \end{cases} 
\end{eqnarray*}
\qed
\end{proof}

\begin{remark} When $\alpha_i = 0$ we can drop the terms with $v_i = v_{i+1}$ in the equation for $[(C_1\times\cdots\times C_n)\cap V_n]$.
\end{remark}

Let us now describe the classes $\iota_*[X_\lambda]$ in $H^*(\PP^5)$ under the Pl\"ucker embedding $\iota:\Gr(1,\PP^3)\hookrightarrow \PP^5$. 
To do so we describe their equations inside $\PP^5$. The degrees of general Schubert varieties were computed by Schubert in \cite{Sch}.

A line $L\in\Gr(1,\PP^3)$ passing through the points $(x_0:x_1:x_2:x_3),(y_0:y_1:y_2:y_3)\in\PP^3$ is uniquely determined by the $2\times 2$ minors of the $2\times 4$ matrix with rows $(x_0:x_1:x_2:x_3)$ and $(y_0:y_1:y_2:y_3)$. Let $p_{i,j}$ denote the minor of the columns $i$ and $j$. The Pl\"ucker embedding associates the vector $(p_{1,2}:p_{1,3}:\ldots:p_{3,4})\in\PP^5$ of the $2\times 2$ minors to each line in $\Gr(1,\PP^3)$.
\begin{itemize}
\item $\Gr(1,\PP^3)$ is defined by the Pl\"ucker relation $p_{1,2}p_{3,4}-p_{1,3}p_{2,4}+p_{2,3}p_{1,4} = 0$. Therefore its degree is $2$ and its codimension $1$, 
so $\iota_*[\Gr(1,\PP^3)]=2z^1$.
\item The condition $L\cap\PP^{(0,1)}\neq \varnothing$ is equivalent to $p_{3,4}=0$. So
$X_1$ is the complete intersection $\{p_{3,4}=0\}\cap\{p_{1,2}p_{3,4}-p_{1,3}p_{2,4}+p_{2,3}p_{1,4}=0\}$. Therefore it has degree $1*2 = 2$ and codimension $2$,
and so $\iota_*[X_1]=2z^2$.
\item Similarly, the condition $L\cap\PP^{(0,1)}\subset\PP^{(0,1,2)}$ is equivalent to $p_{1,4}=p_{2,4}=p_{3,4}=0$. Therefore $X_{1,1}$ is the intersection
\begin{equation*} \{p_{1,2}p_{3,4}-p_{1,3}p_{2,4}+p_{2,3}p_{1,4}=0\}\cap\{p_{1,4}=0\}\cap\{p_{2,4}=0\}\cap\{p_{3,4}=0\}
\end{equation*}
which becomes the complete intersection $\{p_{1,4}=0\}\cap\{p_{2,4}=0\}\cap\{p_{3,4}=0\}$. Therefore it has degree $1*1*1=1$ and codimension $3$,
 and so $\iota_*[X_{1,1}]=z^3$.
\item $X_2$ is the intersection
\begin{equation*} \{p_{1,2}p_{3,4}-p_{1,3}p_{2,4}+p_{2,3}p_{1,4}=0\}\cap\{p_{2,3}=0\}\cap\{p_{2,4}=0\}\cap\{p_{3,4}=0\}
\end{equation*}
which becomes the complete intersection $\{p_{2,3}=0\}\cap\{p_{2,4}=0\}\cap\{p_{3,4}=0\}$. As just above, $\iota_*[X_{2}]=z^3$.
\item $X_{2,1}$ is the intersection
\begin{align*} \{p_{1,2}p_{3,4}-p_{1,3}p_{2,4}+p_{2,3}p_{1,4}=0\}\cap\{p_{1,4}=0\}\cap\{p_{2,3}=0\}\cap\{&p_{2,4}=0\}\\&\cap\{p_{3,4}=0\}
\end{align*}
so the complete intersection $\{p_{1,4}=0\}\cap\{p_{2,3}=0\}\cap\{p_{2,4}=0\}\cap\{p_{3,4}=0\}$. Therefore it has degree $1$ and codimension $4$,
  and so $\iota_*[X_{2,1}]=z^4$.
\item $X_{2,2}$ is the complete intersection $\{p_{1,3}=0\}\cap\{p_{1,4}=0\}\cap\{p_{2,3}=0\}\cap\{p_{2,4}=0\}\cap\{p_{3,4}=0\}$. Therefore it has degree $1$, codimension $5$ and $\iota_*[X_{2,2}]=z^5$.
\end{itemize}
We summarize these computations in Table \ref{table:multideg}.
\begin{table}[h]
\centering
\renewcommand{\arraystretch}{1.5}
\begin{tabular}{| c || c | c |c |c |c |c |c |}\hline
 $X_\lambda$ & $\Gr(1,\PP^3)$ & $X_{1}$& $X_{1,1}$& $X_{2}$& $X_{2,1}$& $X_{2,2}$ \\ \hline
 $\iota_*[X_\lambda]$ & $2z$ & $2z^2$ & $z^3$ & $z^3$ & $z^4$ & $z^5$ \\   
 \hline 
\end{tabular}
\caption{The classes  $\iota_*[X_\lambda]$ in $H^*(\PP^5)$.}
\label{table:multideg}
\end{table}

\begin{remark}
The small case $\Gr(1,\PP^3)$ relevant for this paper is convenient but misleading:
already in $\Gr(1,\PP^4)$ one meets Schubert varieties that are
not complete intersections in the Pl\"ucker embedding, 
making it less straightforward to compute their degrees.
\end{remark}

\begin{proof}[of Theorem \ref{thm:main}]
We compute the multidegrees of $V_n$ and $(C_1\times\cdots\times C_n)\cap V_n$ 
by using Table \ref{table:multideg} to specialize
$$ [X_{2,2}]\mapsto z_i^5, \quad
[X_{2,1}] \mapsto z_i^4,  \quad
 [X_{1,1}],[X_2]\mapsto z_i^3,  \quad
[X_{1}]\mapsto 2z_i^2,  \quad
[\Gr(1,\PP^3)]\mapsto 2z_i $$
in the $i$-th component of Equation (\ref{eq:multivn}).
For $V_n$ we obtain
\begin{eqnarray}
  [V_n] &\mapsto&
  \sum_{0 = v_0 \leq v_1 \leq \ldots \leq v_n = 3
  \atop v_{j+1}-v_j<3} \ 
  \prod_{i=1}^n
  \begin{cases}
    z_i^3 &\text{if }v_i-v_{i-1} = 0 \\
    2z_i^2 &\text{if }v_i-v_{i-1} = 1 \\
    2z_i &\text{if }v_i-v_{i-1} = 2 
  \end{cases}
\end{eqnarray}
Note that given $0 = v_0 \leq v_1 \leq \ldots \leq v_{n-1} \leq v_n = 3$ such that $v_i-v_{i-1}\in\left\{0,1,2\right\}$, there are exactly two different possibilities: either
\begin{enumerate}\item $v_i-v_{i-1}=0$ for all but three indices 
at which $v_i-v_{i-1}=1$, or
\item $v_i-v_{i-1}=0$ for all but two indices $j,k$ 
at which $v_j-v_{j-1}=2$ and $v_k-v_{k-1}=1$.
\end{enumerate}
In the first case, we have a term of the form $8(z_1z_2\cdots z_n)^3z^{-1}_iz^{-1}_jz_k^{-1}$. In the second case we have a term of the form $4(z_1z_2\cdots z_n)^3z_j^{-2}z_k^{-1}$.

Similarly, for $(C_1\times\cdots C_n)\cap V_n$ we obtain
\begin{eqnarray}\label{eq:multimulti}
[(C_1\times\cdots\times C_n)\cap V_n] 
  &=& 
 \hspace{-.4cm} \sum_{0 = v_0 \leq v_1 \leq \ldots \leq v_n = 3
\atop v_{j+1}-v_j<3} \ 
  \prod_{i=1}^n z_i^{5-(v_i-v_{i-1})}
 \hspace{-.15cm} \begin{cases}
    \alpha_i  &\text{if }v_{i} = v_{i-1} \\
    \alpha_i+\beta_i  &\text{if }v_{i}> v_{i-1}
  \end{cases}.
\end{eqnarray}
The formula in the theorem is deduced similarly as in the $[V_n]$ case.
\qed
\end{proof}

\section{$K$-theory}\label{sec:Ktheory}

We conclude this paper by computing the $K$-class of the
concurrent lines variety $V_n$.  Our reference for $K$-polynomials
is \cite[\S 8.5]{MS}, but we include some words of motivation here.

\subsection{A $K$-class example, and the definition}

Consider the diagonal $(\PP^1)_\Delta \subset (\PP^1)^2$, the space
of pairs
$$ \left\{ \left( [x_{0,1},x_{1,1}],[x_{0,2},x_{1,2}] \right)\ :\ 
x_{0,1}/x_{1,1} = x_{0,2}/x_{1,2}
\text{ or really }x_{0,1} x_{1,2} = x_{0,2} x_{1,1}
\right\}. $$
This equation degenerates to $0=x_{0,2} x_{1,1}$ (preserving the homology class), 
which vanishes on the union
$(\PP^1 \times \{0\}) \cup_{(\infty,0)} (\{\infty\} \times \PP^1)$.
Therefore the homology class $[(\PP^1)_\Delta]$ is the sum of the
classes of the two components. (We used the analogous $\PP^2$ 
calculation in the proof of Theorem \ref{thm:cohomVn}.)
\newcommand\calO{{\mathcal O}}
\newcommand\naturals{{\mathbb N}}
\newcommand\integers{{\mathbb Z}}
\newcommand\tensor\otimes

If $A,B$ are subvarieties (of some manifold) of the same dimension,
then the equation $[A\cup B] = [A]+[B]$ we just used in homology
doesn't hold for the ``$K$-classes'' we're about to define: rather,
it obeys a sort of inclusion-exclusion formula
$$ [\calO_{A\cup B}] = [\calO_A] + [\calO_B] - [\calO_{A\cap B}] $$
(the $\calO$s to be defined below). Put another way,
$K$-theory cares about the overcounting of the intersection,
even though its dimension is smaller than that of the components.
In this sense, homology
just sees the top-dimensional part of a $K$-class. 
  Very precisely, there is a filtration on $K(M)$ and an
  isomorphism $gr\ K(M) \otimes \mathbb Q \cong H_*(M;\mathbb Q)$.

We only need to define the $K$-theory of $(\PP^a)^b$, which we do now.
Begin with the abelian group generated by (isomorphism classes of)
finitely generated $\integers^b$-graded modules over the polynomial ring 
in $(a+1)b$ variables $x_{i,j}$, where $x_{i,j}$ acts homogeneously 
with weight $\vec e_j$. We impose two types of relations on this group, the first coming from exact sequences of multigraded modules, and the second
for each $j$, 
coming from modules annihilated by $\langle x_{0,j},\ldots,x_{a,j} \rangle$.
The resulting set of equivalence classes, ``$K$-classes'',
we call $K_0((\PP^a)^b)$.
  Actually Grothendieck derived the name ``$K$-theory''
from the German word ``klasse'', so $K$-class is redundant.

The multigraded modules just described define {\em sheaves} on
$(\PP^a)^b$, i.e. a place for local functions on $(\PP^a)^b$ to
act. (By Liouville's theorem, the only global functions are constant.)
In particular, for $X\subset M$ a closed subscheme (defined by some
multigraded ideal), we know how to multiply local functions on $X$ by
local functions on $M$, giving us a sheaf we call $\calO_X$
with corresponding $K$-class $[\calO_X]$. 
In the $(\PP^1)_\Delta$ example above the two $K$-equivalent modules are
just the quotients by the ideals 
$\langle x_{0,1} x_{1,2} - x_{0,2} x_{1,1} \rangle,
\langle x_{0,2} x_{1,1} \rangle$.

\newcommand\codim{\mathop{\rm codim}}

Defined as above, when $M$ is smooth the abelian group $K_0(M)$ has a 
somewhat non-obvious product, but the only part of it we will need is
$$ [\calO_X] [\calO_Y] = [\calO_{X\cap Y}] \qquad
\text{if the intersection is transverse} $$
for example if $X,Y$ are smooth and $\codim(X\cap Y) = \codim X + \codim Y$.
In $\PP^a$, any two hyperplanes $\PP^{a-1}$ define the same class $H$,
hence $H^{a+1}=0$, 
and it turns out that $K_0(\PP^a) \cong \integers[H]/\langle H^{a+1}\rangle$.

While it will be nice that our methods serve to compute this
finer invariant of $V_n$, our motivation for introducing
$K$-classes is more concrete.
For $X \subseteq \PP^a$ a subscheme, the data of $[\calO_X]$ is 
exactly the same information as the Hilbert polynomial of $X$, 
which one can compute from a Gr\"obner basis for $X$'s ideal.
In particular, a basis $(g_i)$ for $X$'s ideal is Gr\"obner
if {\em and only if} the scheme defined by $\{$the leading
monomials of the $g_i\}$ has the same Hilbert function as $X$,
whose stable behavior is captured by the $K$-class $[X]$.

Let us return to the example of 
$(\PP^1 \times \{0\}) \cup_{(\infty,0)} (\{\infty\} \times \PP^1)$,
whose $K$-class we can now write
as
$$ H \tensor 1 \quad+\quad 1\tensor H \quad-\quad H\tensor H $$
from the inclusion-exclusion formula. We will prefer
to write $H \tensor (1-H) + 1\tensor H$.
We again have, and are using here, a K\"unneth isomorphism
  $K_0((\PP^a)^b) \cong \bigotimes\limits^b \integers[H]/\langle H^{a+1}\rangle$.

\subsection{The $K$-class of $V_n$}

Using a degeneration much like the one we had for $(\PP^1)_\Delta$,
we find the $K$-class $[(\PP^3)_\Delta]_K\in
K((\PP^3)^2)$ to be
\begin{align*}[(\PP^3)_\Delta]_K
&=\sum_{v=0}^3 \left[{\mathcal O}_{\PP^{(0,\cdots,v)}}\right]_K\otimes
\left(\left[{\mathcal O}_{\PP^{(v,\ldots,3)}}\right]_K
-\left[{\mathcal O}_{\PP^{(v+1,\ldots,3)}}\right]_K\right).
\end{align*}
From this formula we can follow arguments similar to those 
of \S \ref{sec:cohom} to compute the $K$-class of $V_n$ in $K((\PP^5)^n)$. 
Letting $H$ denote $[\mathcal O_{\PP^4}]_K \in K(\PP^5)$
the resulting formula is 
\begin{eqnarray*}
  [V_n]_K &=&
  \sum_{0 = v_0 \leq v_1 \leq \ldots \leq v_n = 3
  \atop v_{j+1}-v_j<3} \ 
  \bigotimes_{i=1}^n
  \begin{cases}
    -H^3\otimes(2H-H^2) &\text{if }v_i-v_{i-1} = 0 \\
    (2H^2-H^3)\otimes(2H-3H^2+H^3) &\text{if }v_i-v_{i-1} = 1 \\
    (2H-H^2)\otimes(2H^2-2H^3) &\text{if }v_i-v_{i-1} = 2.
  \end{cases}
\end{eqnarray*}

This is an analogue of the Hilbert polynomial calculation of
\cite[theorem 3.6]{AST}. There are two important differences: theirs
concerns a rational map $\PP^3 \dasharrow (\PP^2)^n$ whereas ours is
about a rational map $\PP^3 \dasharrow Gr(1,\PP^3)^n \to (\PP^5)^n$.
Also, their class in $H^*((\PP^2)^n)$ is multiplicity-free 
in the sense of \cite{Brion}, which is what shows that {\em every}
degeneration of their variety will be reduced (and Cohen-Macaulay),
i.e. that they can have a universal Gr\"obner basis with
squarefree initial terms \cite[\S 2]{AST}.
Our class in $H^*((\PP^5)^n)$ is not multiplicity-free (thanks to 
those coefficients $2$ from Theorem \ref{thm:cohomVn}).

However, Equation (\ref{eq:multivn}) shows that $V_n$ 
{\em is} multiplicity-free in the sense of \cite{Brion} 
when considered as a subvariety of $Gr(1,\PP^3)^n$,
and as such every degeneration of it {\em inside that ambient space}
(not the larger space $(\PP^5)^n$) is reduced and Cohen-Macaulay.

\begin{acknowledgement}
This article was initiated during the Apprenticeship
      Weeks (22 August-2 September 2016), led by Bernd
      Sturmfels, as part of the Combinatorial Algebraic
      Geometry Semester at the Fields Institute for Research in Mathematical Sciences.
      The authors would like to thank their anonymous referees as well as Jenna Rajchgot and Bernd Sturmfels.
      LE was supported by the Fields Institute for Research in Mathematical Sciences.
\end{acknowledgement}

\end{document}